%% file: main.tex
\documentclass[journal, 10pt]{IEEEtran}
\usepackage{amsmath,amssymb,dsfont,graphicx,verbatim}
\graphicspath{{figures/}}
\usepackage{amsmath,amssymb,dsfont,verbatim}
\usepackage{import}
\usepackage{subcaption}
\usepackage{cite}
\usepackage{xcolor}
\usepackage{hyperref}
\usepackage{booktabs} 
\usepackage{tikz}
\usepackage{siunitx}
\usepackage{balance}
\usetikzlibrary{mindmap}

\allowdisplaybreaks[0]

\subimport{custom/}{operators.tex}

\subimport{custom/}{environments.tex}

\title{Tracking Performance of Online Stochastic Learners}
\author{Stefan Vlaski, \IEEEmembership{Member, IEEE}, Elsa Rizk, and Ali H. Sayed, \IEEEmembership{Fellow, IEEE}
\thanks{The authors are with the School of Engineering, \'{E}cole Polytechnique F\'{e}d\'{e}rale de Lausanne. Emails:\{stefan.vlaski, elsa.rizk, ali.sayed\}@epfl.ch.}}
\begin{document}
\maketitle
\begin{abstract}
  The utilization of online stochastic algorithms is popular in large-scale learning settings due to their ability to compute updates on the fly, without the need to store and process data in large batches. When a constant step-size is used, these algorithms also have the ability to adapt to drifts in problem parameters, such as data or model properties, and track the optimal solution with reasonable accuracy.  Building on analogies with the study of adaptive filters, we establish a link between steady-state performance derived under stationarity assumptions and the tracking performance of online learners under random walk models. The link allows us to infer the tracking performance from steady-state expressions directly and almost by inspection.
\end{abstract}
\begin{IEEEkeywords}
Online learning, stochastic learning, tracking performance, non-stationary environment.
\end{IEEEkeywords}
\vspace{-1mm}
\section{Problem Formulation}
\label{sec:intro}
Most online learning algorithms compute an estimate \( \w_i \) at time \( i \) by recursively updating the prior estimate \( \w_{i-1} \) using data \( \x_i \) observed at that same time instant \( i \). We consider in this work a general mapping (i.e., learning rule) of the form:
\begin{equation}\label{eq:def_mapping}
  \w_i = \mathcal{T}(\w_{i-1}; \x_i)
\end{equation}
where \( \mathcal{T}(\cdot, \cdot) \) maps the iterate \( \w_{i-1} \) to \( \w_i \) using the data \( \x_i \). {Throughout this manuscript, we allow for the mapping to be stochastic and time-varying due to the potentially time-varying distribution of the random variable \( \x_i \).} One popular instance of this recursion is the stochastic gradient algorithm~\cite{Polyak87}:
\begin{equation}\label{eq:first_sgd}
  \mathcal{T}(\w_{i-1}; \x_{i}) \triangleq \w_{i-1} - \mu \nabla Q(\w_{i-1}; \x_i)
\end{equation}
which can be used to estimate the minimizer of stochastic risks of the form:
\begin{align}
  \w_i^o \triangleq&\: \argmin_{w\in\mathds{R}^{M}} \mathds{E}_{\x_i} Q(w;\x_i) \label{eq:first_sgd_cost}
\end{align}
where we write \(\w_i^o\), with a subscript \(i\), to allow for the possibility of the minimizer drifting with time due to changes in the distribution of the {streaming} data \(\x_i\). Of course, description~\eqref{eq:def_mapping} captures many more algorithm variations, besides the stochastic gradient algorithm~\eqref{eq:first_sgd}, such as proximal~\cite{Schmidt11,Boyd13}, empirical~\cite{Yuan16}, variance-reduced~\cite{Reddi16, Ge19}, distributed~\cite{Nedic09, Chen12, Chen14multitask, Jaggi14}, and second-order constructions~\cite{Sayed08}. We restrict ourselves in this work to the important class of mappings that satisfy the following mean-square contractive property. We illustrate later by means of examples that several popular learning mappings already satisfy this condition.
\begin{definition}[\textbf{Mean-square contraction}]\label{def:mscontractive}
  We say that a mapping \(\w_i = \mathcal{T}(\w_{i-1}; \x_i) \) is ``mean-square contractive'' around a ``mean-square fixed-point'' \( \w^{\infty}_i \) if for any \(\w_i\) generated by the mapping it holds that:
  \begin{align}\label{eq:general_recursive}
  \E {\left\| \w^{\infty}_{i} - \w_{i} \right\|}^2 \le \gamma_i \E {\left\| \w^{\infty}_{i} - \w_{i-1} \right\|}^2 + \delta_i
  \end{align}
  with \( \gamma_i < 1\). In general, the point \( \w_i^{\infty} \), the rate of contraction \( \gamma_i \), and the additive term \( \delta_i \) will be a function of the distribution of \( \x_i \), and are hence allowed to be time-varying to account for non-stationarity.\hfill\IEEEQED%
\end{definition}
{We refer to the point \( \w_i^{\infty} \) as the ``mean-square fixed-point'' of the mapping \( \mathcal{T}(\w_{i-1}; \x_i) \), since applying \( \mathcal{T}(\cdot; \x_i ) \) at \( \w_{i-1}=\w_i^{\infty} \) yields in light of~\eqref{eq:general_recursive}:
\begin{equation}
  \E {\left\| \w^{\infty}_{i} - \mathcal{T}\left( \w_i^{\infty}; \x_i \right) \right\|}^2 \le \gamma_i \E {\left\| \w^{\infty}_{i} - \w^{\infty}_i \right\|}^2 + \delta_i = \delta_i
\end{equation}
and hence \( \mathcal{T}\left( \w_i^{\infty}; \x_i \right) \approx \w_i^{\infty} \) for small \( \delta_i \) in the mean-square sense.}

If the mapping happens to be deterministic and \( \delta_i = 0 \), we can drop the additive term, as well as the expectation, and recover after taking square-roots:
\begin{equation}
  {\left\| w^{\infty}_{i} - w_{i} \right\|} \le \gamma_i^{\frac{1}{2}} {\left\| w^{\infty}_{i} - w_{i-1} \right\|}
\end{equation}
which corresponds to the traditional definition of a \emph{contractive} mapping~\cite{Kreyszig89}. As we shall show, a number of stochastic algorithms are mean-square contractive, allowing our exposition to cover them all. In the case of the stochastic gradient descent algorithm~\eqref{eq:first_sgd}, the point \( \w_i^{\infty} \) will correspond to the minimizer of~\eqref{eq:first_sgd_cost}, in which case \( \w^o_i \) and \( \w^{\infty}_{i} \) can be used interchangeably. In general, however, such as the decentralized strategies~\eqref{eq:adapt}--\eqref{eq:combine} listed further ahead, we will need to make a subtle distinction.

{In addition to the stochastic nature of the mapping \( \mathcal{T}(\cdot; \x_i) \) resulting from its dependence on the random variable \( \x_i \), we allow for \( \mathcal{T}(\cdot; \x_i) \) to be time-varying due to drifts in the distribution of \( \x_i \), which results in a drift of the fixed-point \( \w_i^{\infty} \) over time (this explains why we are using a subscript \(i\) in \(\w_{i}^{\infty}\)). Relations similar to~\eqref{eq:general_recursive} frequently appear as intermediate results in the performance analysis of stochastic algorithms in \emph{stationary environments}, although stationarity is not necessary for establishing~\eqref{eq:general_recursive}. By establishing a general tracking result for mean-square contractive mappings, {and subsequently appealing to prior results establishing~\eqref{eq:general_recursive},} we can recover known results, and also establish some new results on the tracking performance of stochastic learners for general loss functions.}
\vspace{-2mm}
\subsection{Related Works}
The tracking performance of adaptive filters, focusing primarily on mean-square error designs is fairly well established (see, e.g.,~\cite{Haykin14, Sayed08}). In the decentralized setting, though generally restricted to \emph{deterministic} optimization with exact gradients, the tracking performance of primal and primal-dual algorithms has been studied in~\cite{Ling14, Xi16, Simonetto17}. In the stochastic setting, the tracking performance of the diffusion strategy is established in~\cite{Towfic13}, while the work~\cite{Rizk20} considers a federated learning architecture. The purpose of this work is to establish a unified tracking analysis for the broad class of mean-square contractive mappings, which includes many algorithms as special cases, and will allow us to efficiently recover new tracking results as well.

\section{Tracking Analysis}
\subsection{Non-stationary environments}
We consider a time-varying environment, where the fixed-point \( \w_i^{\infty} \) evolves according to some random-walk model. Such models are prevalent in the study of non-stationary effects.
\begin{assumption}[\textbf{Random Walk}]\label{as:random_walk}
  We assume that the mean-square fixed point of the mapping~\eqref{eq:def_mapping} evolves according to a random walk:
  \begin{equation}\label{eq:random_walk}
    \w_i^{\infty} = \w_{i-1}^{\infty} + \boldsymbol{q}_i
  \end{equation}
  where \( \boldsymbol{q}_i \) is independent of \( \w_{i-1}^{\infty} \). We will allow the random variable \( \boldsymbol{q}_i \) to be non-stationary, with potentially non-zero mean, and only require a global bound on its second-order moment, namely \( \E \|\boldsymbol{q}_i\|^2 \le \xi^2 \).\hfill\IEEEQED%
\end{assumption}
Note that, by allowing \( \boldsymbol{q}_i \) to be non-stationary with non-zero mean, the assumption is more relaxed than typically assumed in the adaptive filtering literature~\cite{Sayed08, Towfic13}. On the other hand, by only imposing a bound on the second-order moment of \( \boldsymbol{q}_i \), rather than on its norm with probability one, condition~\eqref{eq:random_walk} is also more relaxed than in related works on \emph{deterministic} dynamic optimization (e.g.,~\cite{Yuan20}).
Letting \( \widetilde{\w}_{i} \triangleq \w^{\infty}_{i} - \w_{i} \) and using~\eqref{eq:general_recursive}, we have:
\begin{align}\label{eq:drifting_recursive_non_zero}
  \E {\left\| \widetilde{\w}_{i} \right\|}^2 \le&\: \gamma_i \E {\left\|  \widetilde{\w}_{i-1} + \boldsymbol{q}_i \right\|}^2 + \delta_i \notag \\
  \stackrel{(a)}{\le}&\: {\sqrt{\gamma_i}} \E {\left\| \widetilde{\w}_{i-1} \right\|}^2 + \frac{\xi^2}{1 - \sqrt{\gamma_i}} + \delta_i
\end{align}
where in step \( (a) \) we used Jensen's inequality \( \|a+b\|^2 \le \frac{1}{\alpha}\|a\|^2 + \frac{1}{1-\alpha}\|b\|^2 \) for \( 0 < \alpha < 1 \) along with Assumption~\ref{as:random_walk} and \( \gamma_i < 1 \).

If the random variable \( \boldsymbol{q}_i \) happens to be zero-mean and independent of \( \widetilde{\w}_{i-1} \), the inequality can be sharpened by avoiding the use of Jensen's inequality in step \( (a) \) of~\eqref{eq:drifting_recursive_non_zero} and instead appealing to independence of \( \boldsymbol{q}_i \) with \( \widetilde{\w}_{i-1} \) and \( \E \boldsymbol{q}_i = 0 \). This results in:
\begin{align}\label{eq:drifting_recursive}
  \E {\left\| \widetilde{\w}_{i} \right\|}^2 \le \gamma_i \E {\left\| \widetilde{\w}_{i-1} \right\|}^2 +  \xi^2 + \delta_i
\end{align}
In order to continue with the analysis, we assume the following.
\begin{assumption}[\textbf{Global bounds}]\label{as:contraction}
  The rate of contraction \( \gamma_i \) as well as the driving term \( \delta_i \) are bounded from above for all \( i \), i.e., \( \gamma_i \le \Gamma < 1 \) and \( \delta_i \le \Delta \).\hfill\IEEEQED%
\end{assumption}
As we will see in Section~\ref{sec:lms}, Assumption~\ref{as:contraction} generalizes conditions typically imposed in the study of adaptive filters in non-stationary environments. After iterating~\eqref{eq:drifting_recursive_non_zero} and~\eqref{eq:drifting_recursive}, we arrive at the next result.
\begin{theorem}[\textbf{Tracking performance}]\label{TH:TRACKING_PERFORMANCE}
	Suppose \( \mathcal{T}(\cdot; \cdot) \) is a \( (\Gamma, \Delta) \)-mean-square-contractive mapping according to Definition~\ref{def:mscontractive}. Then, we have:
  \begin{equation}\label{eq:performance_non_zero}
    \E \|\widetilde{\w}_i\|^2 \le O\left(\Gamma^{\frac{i}{2}}\right) + \frac{\xi^2}{( 1 - \sqrt{\Gamma} )^2} + \frac{\Delta}{1 - \sqrt{\Gamma}}
  \end{equation}
  In the case when \( \E \boldsymbol{q}_i = 0 \) for all \( i \), we have the tighter relation:
  \begin{equation}\label{eq:performance_zero}
    \E \|\widetilde{\w}_i\|^2 \le O\left(\Gamma^{{i}}\right) + \frac{\xi^2}{1 - {\Gamma}} + \frac{\Delta}{1 - {\Gamma}}
  \end{equation}
\end{theorem}
\begin{IEEEproof}
  The result follows after bounding the quantities appearing in~\eqref{eq:drifting_recursive_non_zero} and~\eqref{eq:drifting_recursive} using Assumption~\ref{as:contraction} and iterating.
\end{IEEEproof}
We note that in steady-state, the terms \(  O(\Gamma^{\frac{i}{2}})  \) and \( O(\Gamma^{{i}}) \) vanish exponentially, and we are left with a drift term proportional to \( \xi^2 \) and a second term proportional \( \Delta \). Furthermore, we note that the non-stationary result~\eqref{eq:performance_zero} can be obtained from the stationary result with \( \xi^2 = 0 \) by merely adding the drift term \( \frac{\xi^2}{1 - {\Gamma}} \).
\vspace{-1mm}
\section{Application to Learning Algorithms}
We now show how Theorem~\ref{TH:TRACKING_PERFORMANCE} can be used to recover the tracking performance of several well-known algorithms under the random walk model~\eqref{eq:random_walk}. We begin by re-deriving and generalizing some known tracking results to illustrate the implications of Assumption~\ref{as:contraction} and verify Theorem~\ref{TH:TRACKING_PERFORMANCE}, and then proceed to derive new tracking results for the multitask diffusion algorithm~\cite{Chen14multitask, Nassif18part1, Nassif20mag}.
\subsection{Least-Mean-Square (LMS) Algorithm}\label{sec:lms}
For illustration purposes, we begin with the least-mean square algorithm, which takes the form:
\begin{equation}\label{eq:lms}
  \mathcal{T}(\w_{i-1}; \boldsymbol{u}_i, \boldsymbol{d}(i)) = \w_{i-1} - \mu \boldsymbol{u}_i \left( \boldsymbol{d}(i) - \boldsymbol{u}_i^{\T}\w_{i-1} \right)
\end{equation}
where the data \( \x_i \triangleq \left\{ \boldsymbol{u}_i, \boldsymbol{d}(i)) \right\} \) arises from the linear model:
\begin{equation}\label{eq:linear_model}
  \boldsymbol{d}(i) = \boldsymbol{u}_i^{\T}\w^{o}_{i} + \boldsymbol{v}({i})
\end{equation}
and \( \boldsymbol{u}_i \in \mathds{R}^M \) denotes an independent sequence of regressors and \( \boldsymbol{v}({i}) \) denotes measurement noise. As is standard in the study of the transient behavior of adaptive filters (see, e.g.,~\cite[Part V]{Sayed08}), we subtract~\eqref{eq:lms} from \( \w_{i}^o \), take squares and expectations to obtain:
\begin{align}
  \E {\left\|\w_i^o - \w_i \right\|}^2 \le \gamma_i  \E \left\| \w_i^o - \w_{i-1} \right\|^2 + \delta_i
\end{align}
with \( R_{u, i} \triangleq \E \boldsymbol{u}_i \boldsymbol{u}_i^{\T} \), \( \gamma_i \triangleq \left\|I - 2 \mu R_{u, i} + \mu^2 \E \boldsymbol{u}_i \boldsymbol{u}_i^{\T}\boldsymbol{u}_i \boldsymbol{u}_i^{\T} \right\| \), \( \sigma_{v, i}^2 \triangleq \E \boldsymbol{v}(i)^2 \) and \( \delta_i \triangleq \mu^2 \mathrm{Tr}\left( R_{u, i} \right) \sigma_{v, i}^2 \). Examination of \( \gamma_i \) and \( \delta_i \) shows that the LMS algorithm~\eqref{eq:lms} satisfies Assumption~\ref{as:contraction} whenever the moments of the regressor \( \boldsymbol{u}_i \) and measurement noise \( \boldsymbol{v}(i) \) are time-invariant (or bounded). This does not restrict the drift of the objective \( \w_i^o \) and the measurement \( \boldsymbol{d}(i) \) which will, of course, be non-stationary as a result. This assumption is also consistent with the modeling conditions typically applied when studying the tracking performance of adaptive filters~\cite[Eq. (20.16)]{Sayed08}. Assuming stationarity of the regressor \( \boldsymbol{u}_i \) and measurement noise \( \boldsymbol{v}(i) \) we find for small step-sizes \( \mu \):
\begin{align}
  \gamma_i \le&\: 1 - 2 \mu \lambda_{\min}\left( R_u \right) + O(\mu^2) \triangleq \Gamma \\
  \delta_i \triangleq&\:  \mu^2 \mathrm{Tr}\left( R_{u, i} \right) \sigma_{v, i}^2 = \mu^2 \mathrm{Tr}\left( R_{u} \right) \sigma_{v}^2 \triangleq \Delta
\end{align}
Hence, we have from~\eqref{eq:performance_zero}:
\begin{align}
  \lim_{i \to \infty} \E {\left\| \widetilde{\w}_{i} \right\|}^2 \le&\:  \frac{ \xi^2}{2 \mu \lambda_{\min}\left( R_u \right) - O(\mu^2)} + \frac{\mu \mathrm{Tr}\left( R_{u} \right) \sigma_{v}^2 }{2 \lambda_{\min}\left( R_u \right) - O(\mu)} \notag \\
  \approx&\:  \frac{ \mu^{-1} \xi^2}{2 \lambda_{\min}\left( R_u \right)} + \frac{\mu \mathrm{Tr}\left( R_{u} \right) \sigma_{v}^2 }{2 \lambda_{\min}\left( R_u \right) }\label{eq:lms_performance}
\end{align}
The result is consistent with~\cite[Lemma 21.1]{Sayed08}, with the factor \( \lambda_{\min}\left( R_u \right) \) appearing in~\eqref{eq:lms_performance} since we are considering here the mean-square deviation of \( \w_i \) around \( \w_i^o \), rather than the excess mean-square error studied in~\cite[Lemma 21.1]{Sayed08}. When the drift term \( \boldsymbol{q}_i \) is no longer zero-mean, we can bound:
\begin{equation}
  \sqrt{\Gamma} = \sqrt{1 - 2 \mu \lambda_{\min}\left( R_u \right) + O(\mu^2)} \le 1 - \mu \lambda_{\min}\left( R_u \right) + O(\mu^2)
\end{equation}
and find from~\eqref{eq:performance_non_zero}:
\begin{align}
  \lim_{i \to \infty} \E {\left\| \widetilde{\w}_{i} \right\|}^2 \le&\:  \frac{ \xi^2}{\left(\mu \lambda_{\min}\left( R_u \right) - O(\mu^2)\right)^2} + \frac{\mu \mathrm{Tr}\left( R_{u} \right) \sigma_{v}^2 }{\lambda_{\min}\left( R_u \right) - O(\mu)} \notag \\
  \approx&\:  \frac{ \mu^{-2} \xi^2}{\lambda_{\min}^2\left( R_u \right)} + \frac{\mu \mathrm{Tr}\left( R_{u} \right) \sigma_{v}^2 }{ \lambda_{\min}\left( R_u \right) }
\end{align}
We observe that the drift penalty incurred in the case when \( \boldsymbol{q}_i \) has non-zero mean is \( O(\mu^{-2}) \), which is significantly larger than in the case where \( \E \boldsymbol{q}_i = 0 \), which is \( O(\mu^{-1}) \). This is to be expected as the cumulative effect of \( \boldsymbol{q}_i \) in the recursive relation~\eqref{eq:random_walk} is no longer equal to zero when \(\E \boldsymbol{q}_i \neq 0 \).

\subsection{Decentralized Stochastic Optimization}
We now consider the problem of general decentralized stochastic optimization. We associate with each agent \( k \) a cost:
\begin{equation}\label{eq:local_costs}
  J_{k, i}(w_k) \triangleq \E Q_{k, i}(w_k; \x_{k, i})
\end{equation}
In this section, we consider the diffusion algorithm for decentralized stochastic optimization~\cite{Chen12, Sayed14}:
\begin{align}
  \boldsymbol{\phi}_{k,i} &= \w_{k,i-1} - \mu {\nabla Q}_{k, i}(\w_{k,i-1}; \x_{k, i})\label{eq:adapt}\\
  \w_{k,i} &= \sum_{\ell=1}^{K} a_{\ell k} \boldsymbol{\phi}_{\ell,i}\label{eq:combine}
\end{align}
for pursuing the minimizer of the aggregate cost:
\begin{equation}\label{eq:pareto}
  \w_i^o \triangleq \argmin_{w}\sum_{k=1}^K p_k J_{k, i}(w)
\end{equation}
where \( \mathrm{col}\{p_k\} \) denotes the right Perron eigenvector associated with the left-stochastic combination matrix \( [A]_{\ell k} = a_{\ell k} \)~\cite{Chen12}. If we collect \( \w_i \triangleq \mathrm{col}\left\{ \w_{k, i} \right\} \) and \( \x_i \triangleq \mathrm{col}\left\{ \x_{k, i} \right\} \), the diffusion recursion~\eqref{eq:adapt}--\eqref{eq:combine} can be viewed as an instance of~\eqref{eq:def_mapping}. Note that by setting the number of agents \( K \) to one we recover ordinary centralized stochastic gradient descent~\eqref{eq:first_sgd}, and as such the results in this section will apply to that case as well. We impose the following standard assumptions on the cost as well as the stochastic gradient approximation~\cite{Sayed14}.
\begin{assumption}[\textbf{Bounded Hessian}]\label{as:hessian}
  Each cost \( J_{k, i}(w) \) is twice-differentiable with bounded Hessian for all \( i \), i.e., \( \nu I \le \nabla^2 J_{k, i}(w) \le \delta I \).\hfill\IEEEQED%
\end{assumption}
Note that this condition ensures that each \( J_{k,i}(\cdot) \) is strongly-convex with Lipschitz gradients and that the respective parameters are \emph{bounded independently of} \( i \). Independence of the bounds on problem parameters over time is common in the study of optimization algorithms in non-stationary and dynamic environments~\cite{Towfic13, Yuan20} and will ensure that Assumption~\ref{as:contraction} is satisfied. We additionally assume that the objectives of the agents do not drift too far apart.
\begin{assumption}[\textbf{Bounded Disagreement}]\label{as:disagreement}
  The distance between each local minimizer is bounded independently of \( i \), i.e.:
  \begin{equation}
    \E \|\w_{k, i}^o - \w_{\ell, i}^o \|^2 \le D^2
  \end{equation}
  for all pairs \( k, \ell\) and times \( i \).\hfill\IEEEQED%
\end{assumption}
 We also make the following common assumption on the quality of the gradient estimate.
\begin{assumption}[\textbf{Gradient noise}]\label{as:gradient_noise}
  Using \( {\nabla Q}_{k, i}(\w_{k,i-1}; \x_{k, i}) \) approximates the true gradient of~\eqref{eq:local_costs} sufficiently well, i.e.:
  \begin{align}
    \E \left\{\nabla Q_{k, i}(w; \x_{k, i})|\boldsymbol{\mathcal{F}}_{i-1}\right\} =&\: \nabla J_{k,i}(w) \\
    \E \left\{\|  \s_{k, i}(w) \|^2 |\boldsymbol{\mathcal{F}}_{i-1}\right\} \le&\: \alpha^2 \| \nabla J_{k, i}(w) \|^2 + \sigma_s^2 \\
      \E \left\{\|  \s_{k, i}(w) \|^2|\boldsymbol{\mathcal{F}}_{i-1}\right\} \le&\: \beta^2 \| w - \w_{k,i}^o \|^2 + \sigma_s^2
  \end{align}
  where \( \boldsymbol{\mathcal{F}}_{i-1} \) denotes the filtration of random variables up to \( i-1\), \( \s_{k, i}(w) \triangleq \nabla Q_{k, i}(w; \x_{k, i}) - \nabla J_{k,i}(w) \), for all \( w \) and some constants \( \alpha^2,  \beta^2, \sigma_s^2 \) independent of \( i \).\hfill\IEEEQED%
\end{assumption}
It has already been established that the diffusion recursion~\eqref{eq:adapt}--\eqref{eq:combine} is a mean-square contractive mapping according to Definition~\ref{def:mscontractive} for some \( \gamma_i \) and \( \delta_i \) in stationary environments~\cite[Eq. (58)]{Chen12}. In order to recover tracking performance through Theorem~\ref{TH:TRACKING_PERFORMANCE}, we need to ensure that the rate of contraction \( \gamma_i \) and driving term \( \delta_i \) can be bounded independent of time \( i \), i.e., that Assumption~\ref{as:contraction} holds under conditions~\ref{as:hessian}--\ref{as:gradient_noise}.
\begin{corollary}[\textbf{Tracking performance of diffusion}]
  The diffusion algorithm~\eqref{eq:adapt}--\eqref{eq:combine} is mean-square contractive around \( \w_{i}^{\infty} \) with \( \Gamma = 1 - 2 \mu \nu + \mu^2 \delta^2 \left( 1 + 4 \alpha^2 \right) \), \( \Delta = \mu^2 N \sigma_s^2 + \mu^2 c_1 \alpha^2 N D^2 \) and
  \begin{align}
    \E \left\| \w_{i}^{\infty} - \mathds{1}\otimes \w_i^o \right\|^2 \le&\: \mu^2 c_2 N \frac{D^2}{1-\lambda_2}\label{eq:diff_bias}
  \end{align}
  where \( \lambda_2 \triangleq \rho\left( A - 1p^{\T} \right) \) denotes the second largest magnitude eigenvalue of the combination matrix \( [A]_{\ell k} = a_{\ell k} \) and \( c_1, c_2 \) denote problem-independent constants. {The quantity \( \w_{i}^{\infty} \) denotes the fixed-point from Definition~\ref{def:mscontractive}, which in light of~\eqref{eq:diff_bias}, is within \( O( \mu^2 )\) of the minimizer of~\eqref{eq:pareto}.} The tracking performance is given by:
  \begin{align}
    \lim_{i\to\infty} \frac{1}{N} \E {\left\|\widetilde{\w}_{k, i}\right\|}^2 \le \mu^{-2} \frac{ 2 \xi^2}{\nu^2} + 2 \mu \frac{\sigma_s^2 + c_1 \alpha^2 D^2}{\nu^2} + \mu^2 \frac{2c_2 D^2}{1-\lambda_2}\label{eq:diffusion_biased}
  \end{align}
  where \( \widetilde{\w}_{k, i} = \w_i^o - \w_{k, i} \). When \( \E \boldsymbol{q}_i = 0 \), we have:
  \begin{align}
    \lim_{i\to\infty} \frac{1}{N} \E {\left\|\widetilde{\w}_{k, i}\right\|}^2 \le \mu^{-1} \frac{ \xi^2}{\nu} + \mu \frac{\sigma_s^2 + c_1 \alpha^2 D^2}{\nu} + \mu^2 \frac{2c_2 D^2}{1-\lambda_2}\label{eq:diffusion_unbiased}
  \end{align}\hfill\IEEEQED%
\end{corollary}
When the gradient approximation \( \nabla Q_{k, i}(w; \x_{k, i}) \) is exact, i.e., \( \alpha^2 = \sigma_s^2 = 0 \), we recover from~\eqref{eq:diffusion_biased} \( \lim_{i\to\infty} \frac{1}{N} \E {\|\w_i^o - \w_{k, i}^2\|}^2 = O(\mu^2) + O(\mu^{-2}) \) which aligns with the result~\cite[Remark 1]{Yuan20}, where \emph{deterministic} dynamic optimization with \emph{exact} gradients is considered. On the other hand, when \( \E \boldsymbol{q}_i = 0 \), we find from~\eqref{eq:diffusion_unbiased} \( \lim_{i\to\infty} \frac{1}{N} \E {\|\w_i^o - \w_{k, i}^2\|}^2 \le O(\mu^{-1}) + O(\mu) + O(\mu^2) \) and recover~\cite[Eq.~(80)]{Towfic13} up to problem-independent factors.

\subsection{Multitask Decentralized Learning}
In this section, we continue to consider a collection of \( K \) agents, each with associated local cost~\eqref{eq:local_costs}. However, instead of pursuing the Pareto solution~\eqref{eq:pareto}, we pursue the multitask problem~\cite{Nassif18part1}:
\begin{equation}\label{eq:multi}
  \w_i^o \triangleq \argmin_{w=\mathrm{col}\{w_k\}} \sum_{k=1}^K J_{k,i}(w_k) + \frac{\eta}{2} w^{\T} \left( L \otimes I \right) w
\end{equation}
where \( L \triangleq \mathrm{diag}\{A \mathds{1}\} - A \) denotes the weighted Laplacian matrix associated with the graph adjacency matrix \( A \). The formulation~\eqref{eq:multi}, in contrast to~\eqref{eq:pareto}, does not force each agent in the network to reach consensus, and instead allows for the independent minimization of \( J_{k,i}(w_k) \) subject to a coupling smoothness regularizer \( \frac{\eta}{2} w^{\T} \left( L \otimes I \right) w \). We refer the reader to~\cite{Nassif18part1, Nassif20mag} for a more detailed motivation for minimizing~\eqref{eq:multi} instead of~\eqref{eq:pareto}, and will focus here on the tracking performance of the resulting algorithm. A solution to~\eqref{eq:multi} can be pursued via the multitask strategy~\cite{Chen14multitask, Nassif18part1}:
\begin{align}
  \boldsymbol{\phi}_{k,i} &= \w_{k,i-1} - \mu {\nabla Q}_{k, i}(\w_{k,i-1}; \x_{k, i})\label{eq:multi_adapt}\\
  \w_{k,i} &= \sum_{\ell=1}^{K} c_{\ell k} \boldsymbol{\phi}_{\ell,i}\label{eq:multi_combine}
\end{align}
where \( c_{\ell k} = 1 - \mu \eta \sum_{\ell = 1}^K a_{\ell k} \) if \( \ell = k \) and \( c_{\ell k} = \mu \eta a_{\ell k} \) otherwise.
Comparing the diffusion strategy~\eqref{eq:adapt}--\eqref{eq:combine} to the multitask strategy~\eqref{eq:multi_adapt}--\eqref{eq:multi_combine} we note a structural similarity with the subtle difference that the combination weights \( c_{\ell k} \) in~\eqref{eq:multi_combine}, in contract to \( a_{\ell k} \) in~\eqref{eq:combine} are \emph{not constant} and depend on the step-size \( \mu \) and regularization parameter \( \eta \). The multitask diffusion strategy~\eqref{eq:multi_adapt}--\eqref{eq:multi_combine} has also been shown to be mean-square contractive~\cite[Eq.~(54)]{Nassif18part1} and hence, we can verify Assumption~\ref{as:contraction} and appeal to Theorem~\ref{TH:TRACKING_PERFORMANCE} to infer its tracking performance.
\begin{corollary}[\textbf{Tracking performance of multitask diffusion}]\label{cor:multi}
  The multitask diffusion algorithm~\eqref{eq:multi_adapt}--\eqref{eq:multi_combine} is mean-square contractive around \( \w_{i}^{\infty} \) with \( \Gamma = 1 - 2 \mu \nu + \mu^2 \left( \delta^2 + 3 \beta^2 \right) \), \( \Delta = \mu^2 N \sigma_s^2 + \mu^2 c_1 3 \beta^2 N D^2 \) and
  \begin{align}
    \E \left\| \w_{i}^{\infty} - \mathds{1}\otimes \w_i^o \right\|^2 \le \mu^2 {\left( \frac{O(\eta^2)}{1+O(\eta^2)} \right)}^2\label{eq:multi_diff_bias}
  \end{align}
  where \( c_1, c_2 \) denote problem-independent constants. {The quantity \( \w_{i}^{\infty} \) denotes the fixed-point from Definition~\ref{def:mscontractive}, which in light of~\eqref{eq:multi_diff_bias}, is within \( O( \mu^2 )\) of the minimizer of~\eqref{eq:multi}.} The tracking performance is hence given by:
  \begin{align}
    \lim_{i\to\infty} \frac{1}{N} \E {\left\|\widetilde{\w}_{k, i}\right\|}^2 \le \mu^{-2} \frac{ 2 \xi^2}{\nu^2} + \mu \frac{2 \sigma_s^2 + 6 c_1 \beta^2 D^2}{\nu^2} + O(\mu^2) \label{eq:multi_diffusion_biased}
  \end{align}
  where \( \widetilde{\w}_{k, i} = \w_i^o - \w_{k, i} \). When \( \E \boldsymbol{q}_i = 0 \), we have:
  \begin{align}
    \lim_{i\to\infty} \frac{1}{N} \E {\left\|\widetilde{\w}_{k, i}\right\|}^2 \le \mu^{-1} \frac{ \xi^2}{\nu} + \mu \frac{\sigma_s^2 + 3 c_1 \beta^2 D^2}{\nu} + O(\mu^2) \label{eq:multi_diffusion_unbiased}
  \end{align}
\end{corollary}
\section{Simulation Results}
\subsection{Tracking Multitask Problems}
\begin{figure}[!t]
	\centering
	\includegraphics[width=\columnwidth]{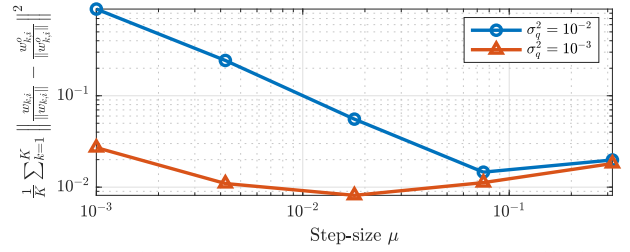}\vspace{-2mm}
	\caption{Tracking performance of the multitask diffusion algorithm~\eqref{eq:multi_adapt}--\eqref{eq:multi_combine} for varying step-sizes \( \mu \) and drift terms \( \sigma_q^2 \).}\label{fig:tracking_performance}
\end{figure}

\vspace{-1mm}
We illustrate the tracking performance of the multitask diffusion strategy~\eqref{eq:multi_adapt}--\eqref{eq:multi_combine} established in Corollary~\ref{cor:multi} in Fig.~\ref{fig:tracking_performance}. We consider a collection of \( K = 20 \) agents connected by a randomly generated graph. Each agent observes feature vectors \( \boldsymbol{h}_{k, i} \) and labels \( \boldsymbol{\gamma}_k(i) \) following a logistic regression model with separating hyperplane \( \w_{k,i}^o \)~\cite[Appendix G]{Sayed14}. The collection of initial hyperplanes \( \{ \w_{k, 0}^o \}_{k=1}^K \) are generated to be smooth over the graph using the procedure of~\cite[Sec. VI]{Nassif18part1} and subsequently subjected to a \emph{common} drift term \( \boldsymbol{q}_i \sim \mathcal{N}\left(0, \sigma_q^2\right) \). Performance is displayed in Fig.~\ref{fig:tracking_performance}. We observe that an optimal step-size choice exists for both drift rates, with smaller \( \xi^2 \) allowing for smaller step-sizes, resulting in smaller effects of the gradient noise and overall better tracking performance. The trends align with the results of Corollary~\ref{cor:multi}.\vspace{-3mm}
\subsection{Illustration of Theorem~\ref{TH:TRACKING_PERFORMANCE} in the Presence of Drift Bias}
\begin{figure}[!t]
	\centering
	\includegraphics[width=\columnwidth]{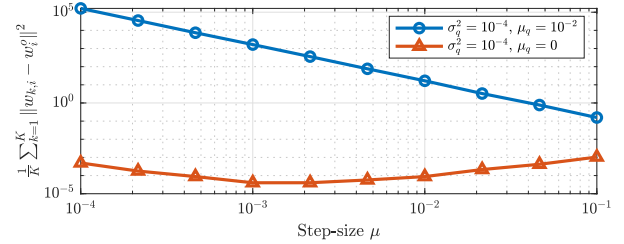}\vspace{-2mm}
	\caption{For \( \E \boldsymbol{q}_i = 0 \), we note reduction in MSD of \( 10  \)\si{dB} per decade for small step-sizes, and a \( 10  \)\si{dB} increase for large step-sizes, which is consistent with \( O(\mu^{-1}) + O(\mu) \) in Eq.~\eqref{eq:diffusion_unbiased}. When \( \mathds{E} \boldsymbol{q}_i \neq 0\), we note a consistent decrease of \( 20  \)\si{dB}, which is consistent with the dominant \( O(\mu^{-2}) \) term in~\eqref{eq:diffusion_biased}.}\label{fig:tracking_biased}
\end{figure}

\vspace{-1mm}
We next verify one of the main conclusions of Theorem~\ref{TH:TRACKING_PERFORMANCE}, namely that the dominant term in the expressions for tracking performance deteriorates from \( O(\mu^{-1}) \) when \( \E \boldsymbol{q}_i = 0 \) (Eq.~\eqref{eq:diffusion_unbiased}) to \( O(\mu^{-2}) \) in the non-zero mean case (Eq.~\eqref{eq:diffusion_biased}). We consider a collection of \( K=5 \) agents observing independent data \( \{ \boldsymbol{u}_{k, i}, \boldsymbol{d}_k(i) \}\) originating from a common linear model \( \w_i^o \in \mathds{R}^3 \) according to~\eqref{eq:linear_model}, subjected to a drift term \( \boldsymbol{q}_i \sim \mathcal{N}\left( \mu_q \mathds{1}, \sigma_q^2 I \right) \). All agents construct local least-squares cost functions \( J_k(w) = \E \|\boldsymbol{d}_k - \boldsymbol{u}_k^{\T} w\|^2 \), and pursue \( \w_i^o \) by means of the resulting diffusion strategy~\eqref{eq:adapt}--\eqref{eq:combine}. The tracking performance in both the zero-mean and biased drift settings for various choices of the step-sizes parameter is displayed in Fig~\ref{fig:tracking_biased}.\vspace{-3mm}

\clearpage
\bibliographystyle{IEEEbib}
{\bibliography{main}}
\balance

\end{document}

%% file: custom/operators.tex
\usepackage{amsthm}

\DeclareMathOperator*{\argmin}{arg\,min}

\DeclareMathOperator{\T}{\mathsf{T}}
\DeclareMathOperator{\E}{\mathds{E}}

\DeclareMathOperator{\w}{\boldsymbol{w}}
\DeclareMathOperator{\x}{\boldsymbol{x}}
\DeclareMathOperator{\s}{\boldsymbol{s}}

%% file: custom/environments.tex
\usepackage{amsthm}

\theoremstyle{plain}

\newtheorem{definition}{Definition}
\newtheorem{assumption}{Assumption}
\newtheorem{theorem}{Theorem}
\newtheorem{corollary}{Corollary}